\documentclass[12pt,a4paper]{article}
\usepackage{amsfonts,amssymb,mathtools,amsxtra,amsthm,array,latexsym,multirow,footnote}
\usepackage{geometry}
\newtheorem{definition}{Definition}%[section]
\newtheorem{proposition}[definition]{Proposition}%[section]

\makesavenoteenv{tabular}
\numberwithin{equation}{section}
\numberwithin{definition}{section}
\newcommand{\Ud}{\mathrm{d}}
\newcommand{\Iff}{if\textcompwordmark f}

\DeclareMathOperator{\rank}{rank}
\DeclareMathOperator{\divz}{div}
\DeclareMathOperator{\sign}{sign}
\DeclareMathOperator{\curl}{curl}

\begin{document}

\title{On polar relative normalizations of ruled surfaces}

\author{{Ioanna-Iris Papadopoulou and Stylianos Stamatakis}
\textbf{\medskip}\\
\emph{Aristotle University of Thessaloniki}\\
\emph{Department of Mathematics}\\
\emph{GR-54124 Thessaloniki, Greece}\\
\emph{e-mail: stamata@math.auth.gr}}
\date{}
\maketitle
\begin{abstract}
\noindent This paper deals with skew ruled surfaces in the Euclidean space $\mathbb{E}^{3}$ which are equipped
 with polar normalizations, that is, relative normalizations such that the relative normal  at each point of the ruled surface lies on the corresponding polar plane. We determine the invariants of a such normalized ruled surface and we study some properties of the Tchebychev vector field and the support vector field of a polar normalization. Furthermore, we study a special polar normalization, the relative image of which degenerates into a curve.
\medskip


\end{abstract}

\section{Introduction} \label{Sec0}
In 1989 F. Manhart introduced the one-parameter family of  relative normalizations $^{(a)} \overline{y}$ of a %relatively normalized
hypersurface with non-vanishing Gaussian curvature $\widetilde{K}$ in the Euclidean space $\mathbb{E}^{n+1}$
% whose relative normalizations
which are characterized by the support functions $^{(a)}q=|\widetilde{K}|^{a}$, $a \in \mathbb{R}$ and called Manhart's normalizations (see \cite{Manhart2}).

G. Stamou and A. Magkos in \cite{Stamou1} and G. Stamou, St. Stamatakis and I. Delivos in \cite{Stamou2} studied
%when a linear combination of the scalar curvature of the relative metric and the relative mean curvature of a ruled surface with non vanishing Gaussian curvature $\widetilde{K}$, which is relatively normalized via the support functions of the form $^{(a)}q=|\widetilde{K}|^{a}$, $a \in \mathbb{R}$ is constant along the generators
ruled surfaces in the  Euclidean space $\mathbb{E}^{3}$ which are equipped with Manhart's normalizations.
 %Some years later G. Stamou, St. Stamatakis and I. Delivos (see \cite{Stamou2}) have studied more properties of a relative normalized ruled surface with non vanishing Gaussian curvature $\widetilde{K}$, which is characterized by the above-mentioned support functions.
Later, S. Stamatakis and I. Kaffas studied in \cite{Stamatakis3} the asymptotic relative normalizations of a ruled surface $\Phi$, that is, relative normalizations such that the relative normals at each point $P$ of $\varPhi$ lie on the corresponding asymptotic plane of $\varPhi$.

Following this idea the authors introduced in \cite{Stamatakis5} three special relative normalizations:
\begin{enumerate}
  \item the \emph{central} normalizations, i.e, relative normalizations such that the relative normals at each point $P$ of $\varPhi$ lie on the corresponding central plane,
  \item the \emph{polar} normalizations, i.e, relative normalizations such that the relative normals at each point $P$ of $\varPhi$ lie on the corresponding polar plane and finally
  \item the \emph{right} normalizations, that is relative normalizations of $\varPhi$ whose relative images $\overline{\varPhi}$ are also ruled surfaces with the additional property that their generators are parallel to those of $\varPhi$. Some of these relative normalizations degenerate into a curve.
\end{enumerate}
The central and the right normalizations were studied thoroughly in \cite{Stamatakis5} and \cite{Stamatakis6}, respectively. In this paper we will study the polar normalizations.

\section{Preliminaries}\label{Sec1}

A brief discussion of some definitions, results and formulae of relative Differential Geometry of surfaces and Differential Geometry of ruled surfaces in the Euclidean space $\mathbb{E}^{3}$ appears in this section. We refer the reader to \cite{Pottmann} and \cite{Schirokow}.

In the three-dimensional Euclidean space $\mathbb{E}^{3}$ let $\varPhi=(U,\overline{x})$ be a ruled $C^{r}$-surface of nonvanishing Gaussian curvature, $r\geq3$, defined by an injective $C^{r} $-immersion $\overline{x}=\overline{x}(u,v)$ on a region $U:=I\times \mathbb{R}$ ($I\subset \mathbb{R}$  open interval)  of $\mathbb{R}^{2}$. We introduce the so-called standard parameters $u\in I, v \in \mathbb{R}$ of $\varPhi$, such that
    \begin{equation}    \label{1}
        \overline{x}(u,v)=\overline{s}(u) + v\,\overline{e}(u),
    \end{equation}%
and
    \begin{equation*}
        \left \vert \overline{e}΄\right \vert =|\overline{e}'|=1,\quad \langle \overline{s}',\overline{e}'\rangle =0,
    \end{equation*}
where the differentiation with respect to $u $ is denoted by a prime and $\langle \,,\rangle$ denotes the standard scalar product in $\mathbb{E}^{3}$..
Here $\varGamma: \overline{s}=\overline{s}(u)$ is the striction curve of $\varPhi$ and the parameter $u$ is the arc length along the spherical curve $\overline{e}=\overline{e}(u)$.

The distribution parameter $\delta(u) :=(\overline{s}',\overline{e},\overline{e}')$, the conical curvature $\kappa(u):=(\overline{e},\overline{e}', \overline{e}'')$ and the function $\lambda(u) :=\cot \sigma$, where $\sigma(u) :=\sphericalangle (\overline{e},\overline{s}')$ is the striction of $\varPhi$ ($-\frac{\pi }{2}<\sigma \leq \frac{\pi }{2}$, $ \sign \sigma = \sign \delta$), are the fundamental invariants of $\varPhi$ and determine uniquely the ruled surface $\varPhi $ up to Euclidean rigid motions.
We also consider the central normal vector $\overline{n}(u):=\overline{e}'$ and the central tangent vector $\overline{z}(u):=\overline{e} \times \overline{n}$. It is known that the vectors of the moving frame $\mathcal{D} : = \{\overline{e}, \overline{n}, \overline{z}\}$ of $\varPhi$ fulfil the following equations \cite[p. 280]{Pottmann}
    \begin{equation} \label{10}
    \overline{e}'=\overline{n},\quad \overline{n}'=-\overline{e}+\kappa \,\overline{z},\quad \overline{z}'=-\kappa \,\overline{n}.
    \end{equation}
Then we have
    \begin{equation}    \label{15}
    \overline{s}'=\delta \,\lambda \, \overline{e}+\delta \, \overline{z}.
    \end{equation}%
We denote partial derivatives of a function (or a vector-valued function) $f$ in the coordinates $u^{1}:=u,\,u^{2}:=v$ by $f_{/i}, f_{/ij}$ etc.
Then from (\ref{1}) and (\ref{15})\ we take
    \begin{equation}\label{20}
    \overline{x}_{/1}=\delta \, \lambda \,\overline{e}+v\,\overline{n}+\delta \,\overline{z},\quad \overline{x}_{/2}=\overline{e},
    \end{equation}
and thus the unit normal vector $\overline{\xi}(u,v) $ to $\varPhi $ is given by
    \begin{equation*}
    \overline{\xi}=\frac{\delta \,\overline{n}-v\,\overline{z}}{w},\quad \text{where}\quad w:=\sqrt{\delta ^{2} + v^{2}}.
    \end{equation*}
Let $I=g_{ij}\Ud u^i \Ud u^j$ and $II=h_{ij}\Ud u^i \Ud u^j$, $i,j=1, 2$ be the first and the second fundamental form of $\varPhi $, respectively, where
\begin{equation}\label{34}
  g_{11}=w^{2} + \delta^{2} \lambda^{2}, \qquad g_{12}=\delta \lambda, \qquad g_{22}=1,
\end{equation}
      \begin{equation}\label{35}
     h_{11}=-\frac{\kappa\, w^{2}+\delta ' \,v-\delta ^{2}\,\lambda }{w},\qquad h_{12}=\frac{\delta }{w}, \qquad  h_{22}=0.
    \end{equation}
The Gaussian curvature $\widetilde{K}(u,v)$ and the mean curvature $\widetilde{H}(u,v)$ of $\varPhi $ are  given by (see \cite{Pottmann})
   \begin{equation}\label{22}
  \widetilde{K}=- \frac{\delta ^{2}}{w^{4}}, \quad \widetilde{H} = -\frac{\kappa w^{2}+\delta ' v+\delta^{2}\lambda}{2w^{3}}.
 \end{equation}
A $C^{s}$-relative normalization of $\varPhi$ is a $C^{s}$-mapping $\overline{y} = \overline{y}(u,v), 1\leq s < r$, defined on $U$, such that
    \begin{equation}\label{45}
    \rank (\{\overline{x}_{/1},\overline{x}_{/2},\overline{y}\})=3,\,\,
    \rank (\{\overline{x}_{/1},\overline{x}_{/2},\overline{y}_{/i}\})=2,\,\,
    i=1,2,\,\,\forall \left(u,v\right) \in U.
    \end{equation}
The pair $\left(\varPhi,\overline {y}\right)$ is called a relatively normalized ruled surface in $\mathbb{R}^{3}$  and the straight line
issuing from a point $P \in \varPhi$ in the direction of $\overline{y}$ is called the relative normal of $\left(\varPhi,\overline {y}\right)$ at $P$. The pair $\overline{\varPhi}=\left(U, \overline{y}\right)$  is called the relative image of $\left(\varPhi,\overline {y}\right)$.

The support function of the relative normalization $\overline{y} $  is defined by $q(u,v):=\langle \overline{\xi},\overline{y}\rangle$ (see~ \cite{Manhart3}).
For $q = 1$, we have $\overline{y} = \overline{\xi}$, that is, the normalization is the Euclidean one.

Due to (\ref{45}), $q $ never vanishes on $U $.
Conversely, when a support function $q$ is given, the relative normalization $\overline{y}$ of the ruled surface $\varPhi$ is uniquely determined and can be expressed in terms of the moving frame $\mathcal{D}$ as follows \cite[p. 179]{Stamatakis3}:
    \begin{equation}\label{70}
    \overline{y}=y_{1}\,\overline{e}+y_{2}\,\overline{n}+y_{3}\,\overline{z},
    \end{equation}
where
    \begin{equation}\label{75}
    y_{1}=-w\frac{\delta q_{/1}+q_{/2}(\kappa \,w^{2}+\delta 'v)}{\delta ^{2}},\quad
    y_{2}=\frac{\delta ^{2}\,q-w^{2}\,v\,q_{/2}}{\delta w},\quad
    y_{3}=-\frac{v\,q+w^{2}\,q_{/2}}{w}.
    \end{equation}
For the coefficients $G_{ij}(u,v)$ of the relative metric $G(u,v)$ of $(\varPhi, \overline{y})$, which is indefinite, we have
 \begin{equation}\label{76}
   G_{ij} = q^{-1} \, h_{ij}.
 \end{equation}
 Then, because of \eqref{35}, the coefficients of the inverse relative metric tensor are
 \begin{equation}\label{77}
   G^{(11)}=0,\quad G^{(12)}=\frac{w\,q}{\delta },\quad G^{(22)}= w\,q \, \frac{\kappa\, w^{2}+\delta 'v-\delta ^{2}\,\lambda }{\delta ^{2}}.
 \end{equation}

For a function (or a vector-valued function) $f$ we denote
 %by $\nabla^{G}\!f$ the first Beltrami differential operator and
 by $\nabla _{i}^{G}f $ the covariant derivatives in the direction of $u^i$, both with respect to the relative metric.
The coefficients $A_{ijk}(u,v)$ of the Darboux tensor are given by
    \begin{equation*}
    A_{ijk}:= q^{-1} \, \langle \overline{\xi},\,\nabla _{k}^{G}\,\nabla _{j}^{G}\,\overline{x}_{/i}\rangle.
    \end{equation*}
Then, by using the relative metric tensor $G_{ij} $ for ``raising and lowering the indices'', the Pick invariant $J(u,v)$ of $(\varPhi,\overline{y})$ is defined by
\begin{equation*}
   J:=\frac{1}{2}A_{ijk}\,A^{ijk}.
\end{equation*}
As we proved in \cite{Stamatakis5} (see equation (2.2)) the Pick invariant is calculated by
 \begin{equation}\label{80}
    J = \frac{3\left(w^2 q_{/2}+v\,q\right)}{2\delta^2 w^3 \,q}\Big\{w^2 \!\left[ \kappa  q v + 2 \delta q_{/1} +q_{/2} \left(\kappa \,w^2 + \delta ' v-\delta^2 \lambda \right)  \right]  -\delta^2 q \left( \lambda v - \delta '  \right) \Big\}.
\end{equation}
The relative shape operator has the coefficients $B_{i}^{j}(u,v)$ given by
    \begin{equation}        \label{90}
    \overline{y}_{/i}=:-B_{i}^{j}\, \overline{x}_{/j}.
    \end{equation}
Then, for the relative curvature $K(u,v)$ and the relative mean curvature $H(u,v)$ of $(\varPhi,\overline{y})$ we have
    \begin{equation}    \label{100}
    K:=\det \left(B_{i}^{j}\right),\quad H:=\frac{B_{1}^{1}+B_{2}^{2}}{2}.
    \end{equation}
We conclude this section by mentioning that, among the surfaces of $\mathbb{E}^3$ with negative Gaussian curvature the ruled surfaces are characterized by the relation
    \begin{equation}\label{136}
    3 H - J -3 S = 0
    \end{equation}
(see \cite{Stamatakis4}), where $S(u,v)$ is the scalar curvature of the relative metric $G$ of such a surface $\varPhi$, which is defined formally as the curvature of the pseudo-Riemannian manifold ($\varPhi,G$).

%The Laplace-normal vector of $\overline{y}$, which was introduced by Heil \cite{Heil} and defined by
%\begin{equation*}
%  \overline{L}=\frac{\vartriangle ^{G} \overline{x}}{2},
%\end{equation*}
%where $\vartriangle ^{G}$ is the second Beltrami-operator with respect to the relative metric $G$, fulfils the relation
%\begin{equation}\label{137}
%  \overline{L}=\overline{T}+\overline{y}.
%\end{equation}

\section{Polar normalizations}\label{Sec2}
We concentrate now on the main topic of this paper, namely the polar normalizations of a skew ruled surface $\varPhi$, i.e., relative normalizations  such that the relative normal at each point $P$ of $\varPhi$ lies on the corresponding polar plane $\{ P;\overline{n}, \overline{z} \}$. In \cite{Stamatakis5} it was shown that the support function of $\overline{y}$ is of the form
\begin{equation}\label{110}
  q=f(V),
\end{equation}
where $f(V)$ is an arbitrary $C^{2}$-function of
\begin{equation}\label{115}
  V= \arctan \frac{v}{\delta}- \!\! \int \! \kappa \Ud u.
\end{equation}
By means of \eqref{70}, \eqref{75}, \eqref{110} and \eqref{115} we deduce that the arising relative normalization, i.e., the polar normalization of the given ruled surface $\varPhi$ is
\begin{equation}\label{116}
\overline{y}=\frac{\delta q- \dot{q} v}{w} \overline{n}- \frac{q v+ \delta \dot{q}}{w}\overline{z},
\end{equation}
where the dot denotes the differentiation with respect to $V$.
Then, from \eqref{10}, \eqref{20}, \eqref{90} and \eqref{116}, we take the coefficients $B_i^j$ of the relative shape operator of a polar normalization:
\begin{align*}
  B_{1}^{1} &= - \frac{(\kappa w^{2}+\delta ' v)(q+\ddot{q})}{w^{3}},\\
  B_{1}^{2} &= \frac{1}{w^{3}} \Big \{ - \dot{q} v^{3} -\delta^{2} \dot{q} v + \delta^{3}\left[ q\left(\kappa \lambda +1\right) + \kappa \lambda \ddot{q}\right] + \delta v \left[q\left(\kappa \lambda v+ v+ \delta ' \lambda\right)+ \lambda \ddot{q}\left(\kappa v + \delta' \right)\right] \Big \}, \\
   B_{2}^{1}& = \frac{\delta (q+\ddot{q})}{w^{3}},\\
  B_{2}^{2} &= - \frac{\delta^{2} \lambda (q+\ddot{q})}{w^{3}}.
\end{align*}
Hence, by using \eqref{100} and (\ref{22}b), we obtain the relative curvature $K$ and the relative mean curvature $H$:
\begin{equation}\label{120}
  K=-\delta \, \frac{\left(\delta q- \dot{q} v \right) \left(q+ \ddot{q} \right)}{w^{4}}, \quad H=\widetilde{H}\left(q+ \ddot{q}\right).
\end{equation}
From (\ref{120}a) we deduce that the relative curvature $K$ of a polar normalization vanishes identically \Iff{}
 \begin{equation*}
   \delta q- \dot{q} v=0 \quad \text{or} \quad  q+ \ddot{q}=0,
 \end{equation*}
 or, equivalently, \Iff{}
 \begin{equation*}
   q=c e^{\frac{\delta V}{v}}, c\in\mathbb{R^{*}} \quad  \text{or} \quad q= c_{1} \cos V + c_{2} \sin V, \, c_{1}, c_{2} \in\mathbb{R},\,  c_{1}^{2}+c_{2}^{2} \neq 0.
 \end{equation*}
 We reject the first support function since it leads to a non polar normalization. Thus we have the following

\medskip
 \begin{proposition}
   Let $\varPhi \subset E^{3}$ be a polar normalized ruled surface. The relative curvature $K$ of $(\varPhi,\overline{y})$ vanishes identically \Iff{} the support function is of the form
  \begin{equation*}
  q= c_{1} \cos V + c_{2} \sin V, \quad c_{1}, c_{2} \in\mathbb{R},\quad  c_{1}^{2}+c_{2}^{2} \neq 0.
\end{equation*}
 \end{proposition}

By taking (\ref{22}b) and (\ref{120}b) into consideration we arrive at

 \medskip
 \begin{proposition}
   Let $\varPhi \subset E^{3}$ be a polar normalized ruled surface. $(\varPhi,\overline{y})$ is relatively minimal $(H=0)$ \Iff{} one of the following holds true \\
    \textup{(a)} the support function is of the form
      \begin{equation*}
      q= c_{1} \cos V + c_{2} \sin V, \quad c_{1}, c_{2} \in\mathbb{R},\quad  c_{1}^{2}+c_{2}^{2} \neq 0, \\
        \end{equation*}
     \textup{(b)} $(\varPhi,\overline{y})$ is a polar normalized right helicoid ($\delta=c\in\mathbb{R^{*}}$ and $\kappa=\lambda=0$).
 \end{proposition}

We notice that both the relative curvature $K$ and the relative mean curvature $H$ vanish identically \Iff{} the support function is of the form
\begin{equation}\label{122}
  q= c_{1} \cos V + c_{2} \sin V, \quad c_{1}, c_{2} \in\mathbb{R},\quad  c_{1}^{2}+c_{2}^{2} \neq 0.
\end{equation}

By using (\ref{22}b) and \eqref{80} we find the Pick invariant
\begin{equation}\label{125}
  J = \left( q v + \delta \, \dot{q} \right) \left( \frac{J_{EUK}}{v} + \frac{3 \widetilde{H} \,\dot{q}}{\delta \, q}\right),
\end{equation}
where
\begin{equation}\label{126}
J_{EUK}=3 v \, \frac{\kappa v^{3} + \delta^{2}\left(\kappa-\lambda\right)v+\delta^{2} \delta'}{2\delta^{2} w^{3}}
\end{equation}
is the Pick invariant of the Euclidean normalization. The Pick invariant vanishes identically \Iff{}
 \begin{equation*}
   q v + \delta \, \dot{q}=0 \quad \text{or} \quad \frac{J_{EUK}}{v} + \frac{3 \widetilde{H} \,\dot{q}}{\delta \, q}=0,
 \end{equation*}
  or, equivalently, \Iff{}
  \begin{itemize}
    \item the support function is of the form
    \begin{equation*}
      q=c_{1} e^{\frac{-Vv}{\delta}}, \,c_{1}\in\mathbb{R^{*}},\, \text{or}
    \end{equation*}
    \item $\varPhi$ is not a right helicoid and the support function is of the form
    \begin{equation*}
      q= c_{2} e^{\frac{V\left[\kappa v^{3}+ \delta^{2}\left(\kappa-\lambda\right)v+\delta^{2} \delta\,'\right]}{\delta\left[\kappa v^{2}+\delta\,'v+\delta^{2}\left(\kappa+\lambda\right)\right]}}, \, \,c_{2}\in\mathbb{R^{*}},\, \text{or}
    \end{equation*}
    \item $\varPhi$ is a right helicoid.
  \end{itemize}

We reject the two support functions since they are not polar. Hence, we deduce

\medskip
 \begin{proposition}
   Let $\varPhi \subset E^{3}$ be a polar normalized ruled surface. The Pick invariant $J$ of $(\varPhi,\overline{y})$ vanishes identically \Iff{} $\varPhi$ is a right helicoid.
 \end{proposition}

From \eqref{136}, (\ref{120}b), \eqref{125} and \eqref{126} we evaluate the scalar curvature of the relative metric
\begin{equation*}
\begin{split}
    S& =\frac{1}{2 \delta^{2} w^{3} q} \Big\{ -q^{2}\left\{\kappa w^{4}+\delta^{2}\left[\left(-v^{2}+\delta^{2}\right)\lambda+2\delta'v\right]\right\}+\delta^{2}\left(\kappa w^{2}+\delta^{2}\lambda+\delta' v\right)\dot{q}^{2} \\
     & + \delta q \big \{\left[2 \delta^{2}\lambda v+\left(v^{2}-\delta^{2}\right)\delta'\right]\dot{q}-\delta\left(\kappa w^{2}+\delta'v+\delta^{2}\lambda\right)\ddot{q}\big\}\Big\}.
\end{split}
 \end{equation*}

\section{The Tchebychev vector field and the support vector field of a polar normalization}\label{Sec3}
In \cite{Stamatakis3} it was shown that the coordinate functions of the Tchebychev vector $\overline{T}(u,v)$ of $(\varPhi,\overline{y})$, which is defined by
\begin{equation*}
 \overline{T}:=T^{m}\, \overline{x}_{/m},\quad \text{where\quad }T^{m}:=\frac{1}{2}A_{i}^{im},
 \end{equation*}
are given by
 \begin{equation*}
 T^{1}=\frac{w^{2}q_{/2} + v\,q}{\delta \,w},\,\, T^{2}=\frac{2\delta\,w^{2}q_{/1}+\delta'q\,(\delta^{2}-v^{2})}{2\delta^{2}\, w}+\frac{T^{1}(\kappa w^{2}+\delta'v-\delta^{2}\lambda)}{\delta}.
\end{equation*}
Hence, by using \eqref{110} and \eqref{115}, we deduce that the coordinate functions of the Tchebychev vector of a polar normalization are
\begin{equation}\label{200}
T^{1}=\frac{q v+\delta\dot{q}}{\delta w}, \quad T^{2}=\frac{q\left(2 \kappa v w^{2}-2 \delta^{2} \lambda v+ \delta' w^{2}\right)-2 \delta^{3}\lambda\dot{q}}{2 \delta^{2}w}.
\end{equation}
The divergence $\divz^{I}\overline{T}$ of $\overline{T}$ with respect to the first fundamental form $I$ of $\varPhi$, which initially reads (see \cite{Stamatakis3})
\begin{equation}\label{205}
\divz^{I}\overline{T}=\frac{\left(w \, T^{i} \right)_{/i}}{w},
\end{equation}
becomes, on account of \eqref{200},
\begin{equation*}
\begin{split}
   \divz^{I}\overline{T} & = \frac{1}{2 \delta^{2} w^{3}} \bigg\{ 2w^{2}q \Big[ \left( 3v^{2} + \delta^{2} \right) \kappa-\delta^{2}\lambda \Big] \\
     & + \delta \Big\{ \Big[ -\delta' v^{2}+\delta^{2}\left(-2\lambda v+ \delta'\right)\Big]\dot{q}-2\delta\left(\kappa w^{2}+\delta^{2}\lambda+\delta'v\right)\ddot{q}\Big\} \bigg\}.
\end{split}
\end{equation*}
The rotation $\curl^{I}\overline{T}$ of $\overline{T}$ with respect to the first fundamental form $I$ of $\varPhi$, which initially reads (see \cite{Stamatakis3})
\begin{equation}\label{207}
\curl^{I}\overline{T}=\frac{\left(g_{12}T^{1}+g_{22}T^{2}\right)_{/1}-\left(g_{11}T^{1}+g_{12}T^{2}\right)_{/2}}{w},
\end{equation}
becomes, by taking \eqref{34} and \eqref{200} into consideration,
\begin{equation*}
  \begin{split}
     \curl^{I}\overline{T} & = - \frac{1}{2 \delta^{3}w^{2}} \bigg \{2 \delta'q v^{2}\left(2 \kappa v+\delta'\right)+\delta^{2} q\Big[ 4 \left(\kappa \lambda +1\right)v^{2}+\delta'\left(2\kappa+\lambda\right)v+\delta'^{2}\Big] \\
       & +  \Big \{ \dot{q} \Big[ 4v+\left(\kappa+\lambda\right)\left(2\kappa v+\delta'\right)\Big]-q\left(2\kappa'v+\delta''\right)\Big \}\\
       & + \delta v \Big[ 2 \kappa^{2}\dot{q} \,v^{2}+3\kappa\delta'\dot{q}\,v+\delta'^{2}\dot{q}-q v \left(2\kappa' v+ \delta''\right)\Big]+ 2\delta^{4}\Big[ q\left(\kappa \lambda +1\right)+\ddot{q}\Big] \bigg \}.
  \end{split}
\end{equation*}
%The divergence $\divz^{G}\overline{T}$ of $\overline{T}$ with respect to the relative metric of $\varPhi$, which initially reads
%\begin{equation}\label{208}
%\divz^{G}\overline{T}=\frac{\left(G T^{i} \right)_{/i}}{|G|^{1/2}},
%\end{equation}
%becomes, by means of \eqref{76} and \eqref{200},
%\begin{equation*}
%  \begin{split}
%     \divz^{G}\overline{T} &= \frac{1}{\delta^{2}w^{3}q}\Big \{ q^{2} \big \{ \kappa w^{4}+ \delta^{2} \left[\left(v^{2}-\delta^{2}\right) \lambda - 2 \delta' v \right] \big \} +\delta^{2} \dot{q}^{2} \left(\kappa w^{2} + \delta ' v+ \delta^{2} \lambda\right)\\
%       & + \delta q \big \{\dot{q} \left[2 \delta^{2}\lambda v+ \delta' \left(v^{2}-\delta^{2}\right)\right]-\delta\ddot{q} \left(\kappa w^{2}+\delta' v+ \delta^{2}\lambda\right) \big \} \Big\}.
%  \end{split}
%\end{equation*}
%The rotation $\curl^{G}\overline{T}$ of $\overline{T}$ with respect to the relative metric of $\varPhi$, which initially reads
%\begin{equation}\label{209}
%\curl^{G}\overline{T}=\frac{\left(G_{12}T^{1}+G_{22}T^{2}\right)_{/1}-\left(G_{11}T^{1}+G_{12}T^{2}\right)_{/2}}{|G|^{1/2}},
%\end{equation}
%becomes, by using \eqref{35}, \eqref{76}, \eqref{110}, \eqref{115} and \eqref{200},
%\begin{equation*}
%  \curl^{G}\overline{T}=0,
%\end{equation*}

Analogously we calculate the divergence and the rotation of $\overline{T}$ with respect to the relative metric of $\varPhi$:
\begin{align*}
  \begin{split}
     \divz^{G}\overline{T} &= \frac{1}{\delta^{2}w^{3}q}\Big \{ q^{2} \big \{ \kappa w^{4}+ \delta^{2} \left[\left(v^{2}-\delta^{2}\right) \lambda - 2 \delta' v \right] \big \} +\delta^{2} \dot{q}^{2} \left(\kappa w^{2} + \delta ' v+ \delta^{2} \lambda\right)\\
       & + \delta q \big \{\dot{q} \left[2 \delta^{2}\lambda v+ \delta' \left(v^{2}-\delta^{2}\right)\right]-\delta\ddot{q} \left(\kappa w^{2}+\delta' v+ \delta^{2}\lambda\right) \big \} \Big\},
  \end{split}\\
  \curl^{G}\overline{T}& = 0.
\end{align*}
Last relation agrees with
\begin{equation*}
  \overline{T}=\nabla^{G} \left( \ln \frac{q}{q_{AFF}},\overline{x}\right) \quad \text{where} \quad q_{AFF}= | \widetilde{K}|^{1/4}
\end{equation*}
(see \cite{Stamatakis4}).
So, we have

\medskip
 \begin{proposition}
   Let $\varPhi \subset E^{3}$ be a polar normalized ruled surface. The rotation of the Tchebychev vector field with respect to the relative metric of $\varPhi$ vanishes identically and its potential is given by
   \begin{equation*}
  \tau \left(u, v\right)= \ln \frac{w q}{|\delta|^{1/2}}+c, \quad c\in\mathbb{R}.
\end{equation*}
 \end{proposition}

Now let
\begin{equation}\label{211}
  \overline{Q}:= \frac{1}{4} \bigtriangledown^{G}\!\Big(\frac{1}{q},\overline{x}\Big)
\end{equation}
be the support vector $\overline{Q}(u,v)$ of $( \varPhi, \overline{y})$, which is introduced in  \cite{Stamatakis3}. By taking \eqref{77}, \eqref{110} and \eqref{115} into consideration we find that the coordinate functions of the support vector field of a polar normalization are
\begin{equation}\label{210}
  Q^{1}=-\frac{\dot{q}}{4wq}, \quad Q^{2}=\frac{\delta \lambda \dot{q}}{4wq}.
\end{equation}
By means of (\ref{22}b), \eqref{205} and \eqref{210}, we find the divergence $\divz^{I}\overline{Q}$ of $\overline{Q}$ with respect to the first fundamental form $I$ of $\varPhi$
\begin{equation*}
  \divz^{I}\overline{Q}= \widetilde{H}\frac{\dot{q}^{2}-q \ddot{q}}{2q^{2}}.
\end{equation*}
Hence, we derive

 \medskip
 \begin{proposition}
   Let $\varPhi \subset E^{3}$ be a polar normalized ruled surface. The support vector field is incompressible with respect to the first fundamental form of $\varPhi$ $(\divz^{I}\overline{Q}=0)$ \Iff{} \\
    \textup{(a)} the support function is of the form
     \begin{equation*}
  q=c_{2}e^{c_{1}V}, \quad c_{1} \in\mathbb{R}, \quad c_{2} \in\mathbb{R^{*}}, \, \, \text{or} \\
\end{equation*}
    \textup{(b)} $\varPhi$ is a right helicoid.
 \end{proposition}

By taking \eqref{34}, \eqref{207} and \eqref{210} into account we deduce that the rotation $\curl^{I}\overline{Q}$ of $\overline{Q}$ with respect to the first fundamental form $I$ of $\varPhi$ is
\begin{equation*}
  \curl^{I}\overline{Q}=\frac{-\delta \dot{q}^{2}+q\left(\dot{q}v+\delta\ddot{q}\right)}{4w^{2}q}.
\end{equation*}
By taking \eqref{76}, \eqref{205} and \eqref{210} into consideration we find the divergence $\divz^{G}\overline{Q}$ of $\overline{Q}$ with respect to the relative metric of $\varPhi$
\begin{equation*}
\begin{split}
   \divz^{G}\overline{Q} &  = \frac{1}{4 \delta w^{3} q^{2}} \Big\{ \dot{q}\big \{q \left[-\delta' v^{2}+\delta^{2}\left(-2\lambda v+ \delta'\right)\right]-2 \delta \dot{q} \left(\kappa w^{2}+\delta'v+\delta^{2}\lambda\right) \big \} \\
     & + \delta q \ddot{q} \left(\kappa w^{2}+ \delta^{2}\lambda+ \delta' v\right) \Big \}.
\end{split}
 \end{equation*}
By using \eqref{35}, \eqref{76}, \eqref{110}, \eqref{115}, \eqref{207} and \eqref{210} we have the rotation $\curl^{G}\overline{Q}$ of $\overline{Q}$ with respect to the relative metric of $\varPhi$
\begin{equation*}
\curl^{G}\overline{Q}=0,
\end{equation*}
which agrees with the relation \eqref{211}. Thus, we have

\medskip
 \begin{proposition}
   Let $\varPhi \subset E^{3}$ be a polar normalized ruled surface. The rotation of the support vector field with respect to the relative metric of $\varPhi$ vanishes identically and its potential is given by
   \begin{equation*}
  \tau \left(u, v\right)= \frac{1}{4 q}+c, \quad c\in\mathbb{R}.
\end{equation*}
 \end{proposition}

\section{A special polar normalization}\label{Sec4}
In this section we will study the support function of the form \eqref{122}, which arises when the relative curvature $K$ or the relative mean curvature $H$ vanishes identically (see Sec. \ref{Sec2}). By using \eqref{116} the corresponding relative normalization takes the form
\begin{equation*}
  \overline{y}=\left[c_{1} \cos\left( \, \, \! \! \int \! \kappa \Ud u \right)-c_{2} \sin \left( \, \, \! \! \int \! \kappa \Ud u \right) \right] \overline{n} - \left[c_{2} \cos\left( \, \, \! \! \int \! \kappa \Ud u \right)+c_{1} \sin \left( \, \, \! \! \int \! \kappa \Ud u \right) \right] \overline{z},
\end{equation*}
i.e., the relative normalization degenerates into a curve $\varGamma^{*}$ with curvature
\begin{equation*}
  \kappa^{*}=\frac{1}{|c_{1} \cos\left( \, \, \! \! \int \! \kappa \Ud u \right)-c_{2} \sin \left( \, \, \! \! \int \! \kappa \Ud u \right)|}
\end{equation*}
and torsion
\begin{equation*}
  \sigma^{*}=\frac{-\kappa}{c_{1} \cos\left( \, \, \! \! \int \! \kappa \Ud u \right)-c_{2} \sin \left( \, \, \! \! \int \! \kappa \Ud u \right)}.
\end{equation*}
Since
\begin{equation*}
  \frac{\kappa^{*}}{\sigma^{*}}=\pm \frac{1}{\kappa},
\end{equation*}
 we deduce that \emph{$\overline{y}$ is a curve of constant slope \Iff{} $\varPhi$ is a ruled surface of constant slope}.

By means of \eqref{125} and \eqref{126} we find the Pick invariant of this normalization:
\begin{equation*}
  \begin{split}
     J & = \frac{3\left[c_{2} \cos\left( \, \, \! \! \int \! \kappa \Ud u \right)+c_{1} \sin \left( \, \, \! \! \int \! \kappa \Ud u \right) \right]}{2 \delta^{2} w \left(c_{1} \cos V +c_{2} \sin V \right)} \Big \{ \cos \left( \, \, \! \! \int \! \kappa \Ud u \right) \Big [\kappa \left(c_{2} v^{2} + 2 c_{1} \delta v - c_{2} \delta^{2}\right) \\
       & + \delta \left(-c_{2} \delta \lambda + c_{1} \delta' \right) \Big ]+ \sin \left( \, \, \! \! \int \! \kappa \Ud u \right) \Big
       [ \kappa \left(c_{1} v^{2} - 2 c_{2} \delta v - c_{1} \delta^{2} \right)- \delta \left(c_{1} \delta \lambda + c_{2} \delta '\right) \Big ] \Big \}.
  \end{split}
\end{equation*}
Then by using \eqref{20} and \eqref{200} we deduce the Tchebychev vector
\begin{equation*}
\begin{split}
   \overline{T}& = \frac{w}{2 \delta^{2}} \left(c_{1} \cos V + c_{2} \sin V \right) \left(2 \kappa v + \delta' \right) \overline{e}+ \frac{v}{\delta}\left[c_{2} \cos\left( \, \, \! \! \int \! \kappa \Ud u \right)+c_{1} \sin \left( \, \, \! \! \int \! \kappa \Ud u \right) \right]\overline{n} \\
     & + \left[c_{2} \cos\left( \, \, \! \! \int \! \kappa \Ud u \right)+c_{1} \sin \left( \, \, \! \! \int \! \kappa \Ud u \right) \right] \overline{z}.
\end{split}
\end{equation*}
Finally, by taking \eqref{20} and \eqref{210} into consideration we derive the support vector
\begin{equation*}
  \overline{Q}= \frac{c_{1} \sin V - c_{2} \cos V}{4w \left(c_{1} \cos V + c_{2} \sin V\right)} \left(v \overline{n}+ \delta \overline{z}\right).
\end{equation*}
%Inserting \eqref{400} and \eqref{405} in \eqref{137} we derive the Laplace-normal vector
%\begin{equation*}
%\begin{split}
%    \overline{L}&  = \frac{w}{2 \delta^{2}} \left(c_{1} \cos V + c_{2} \sin V \right) \left(2 \kappa v + \delta' \right) \overline{e}  \\
%     & + \frac{\cos\left( \, \, \! \! \int \! \kappa \Ud u \right) \left(c_{2}v+c_{1}\delta\right)+\sin \left( \, \, \! \! \int \! \kappa \Ud u \right)\left(c_{1}v-c_{2}\delta\right)}{\delta}\overline{n}.
%\end{split}
%\end{equation*}

\end{document}